\documentclass[11pt]{article}
\usepackage{graphicx}
\usepackage{amsmath,amssymb,multirow,longtable,rotating,array,booktabs,float,bm,amsthm,diagbox}
\usepackage{pdflscape}
\usepackage{multirow,booktabs,makecell}
\usepackage{epstopdf}
\usepackage[colorlinks,linkcolor=red,anchorcolor=blue,citecolor=green]{hyperref}
\allowdisplaybreaks
\newtheorem{theorem}{Theorem}

\catcode`@=11 \@addtoreset{equation}{section} \catcode`@=12

\hypersetup{
    colorlinks=true,
    linkcolor=blue,
    urlcolor=blue,
    citecolor=blue,
}
\usepackage{cases}
\begin{document}
	\title{Bell-INGARCH Model.
		\footnotetext{$\dag$Corresponding author. School of Mathematics and Statistics, Jiangsu Normal University, Xuzhou, 221116, China, E-mail:  qianly92@163.com}}
	\author{  Ying Wang$^1$, Shuang Chen$^2$ and Lianyong Qian$^{2\dag}$\\[0cm]
		{\small\it $^1$School of Mathematics and Physics, China University of Geosciences, Wuhan, 100083, China}\\
		{\small\it $^2$School of Mathematics and Statistics, Jiangsu Normal University, Xuzhou, 221116, China}}
	\date{}
	\maketitle
	\begin{center}
		\begin{minipage}{14.5truecm}
			{\bf Abstract}
 Integer-valued time series exist widely in economics, finance, biology, computer science, medicine, insurance, and many other fields. In recent years, many types of models have been proposed to model integer-valued time series data, in which the integer autoregressive model and integer-valued GARCH model are the most representative. Although there have been many results of integer-valued time series data, the parameters of integer-valued time series model structure are more complicated. This paper is dedicated to proposing a new simple integer-valued GARCH model. First, the Bell integer-valued GARCH model is given based on Bell distribution. Then, the conditional maximum likelihood  estimation method is used to obtain the estimators of parameters. Later, numerical simulations confirm the finite sample properties of the estimation of unknown parameters. Finally, the  model is applied in the two real examples. Compared with the existing models, the proposed model is more simple and applicable.

   
\end{minipage}
\end{center}
%%%%%%%%%%%%%%%%%%%%%%%%%%%%%%%%%%%%%%%%%%%%%%%%%%%%%%%%%%%%%%%%%%%%%%%%%%%%%%%%%%%%%%%%%%%%%%%%%%%%%%%%%%%%%%

\section{Introduction }
Integer-valued time series data are widely applied in various fields such as economics, finance, biology, computer science, medicine, insurance, and more. For example, the weekly count of new virus infections in a particular region during a period, the monthly crime rates in a city or region over several years, the number of stock trades in five-minute intervals within a day, and so on. Integer-valued time series data have attracted increasing attention from scholars and have become a significant research area in modern statistics.

Traditional time series models typically deal with continuous data, and these models fail to capture specific characteristics of integer-valued time series, such as asymmetric marginal distributions and zero inflation. Common issues include an inability to interpret data generation mechanisms, substantial biases when approximating with continuous data, and often non- integer forecast values. As a result, many researchers have introduced new models. Currently, two types of models are primarily used to model integer-valued time series data. One type is the integer-valued autoregressive model (INAR) based on the sparse operators, and the other type is the integer-valued GARCH (INGARCH).

For count data, researchers are initially devoted to the Poisson distribution. As research delved more profound, it was found that the Poisson distribution was only suitable for handling equi-disperesed count data. Many real examples, however, exhibit overdispersion, meaning their variance is greater than their mean. Consequently,  alternative distributions are considered as the consumptions for the INAR and INGARCH models. For instance, Jazi (2005) introduced an INAR model with zero-inflated Poisson innovations. Qian (2020) proposed an INAR model featuring innovations with a generalized Poisson-inverse Gaussian distribution. Ferland (2006) proposed the Poisson INGARCH model (PINGARCH), which defined as follows:
\begin{equation}{\label{e1.1}}
X_{t} | \mathcal{F}_{t-1}\sim\mathcal{P}(\mathcal{\lambda}_{t}),                      
   \mathcal{\lambda}_{t}=\alpha_0+\alpha_1 X_{t-1}+\beta_1 \mathcal{\lambda}_{t-1},
\end{equation}
where $\alpha_0>0$, $\alpha_1\geq 0$, $\beta_1\geq 0$ and $\mathcal{F}_{t-1}$ is the $\sigma$-fields generated by $\{X_s: s\le t-1\}$.
Zhu (2011) introduced the negative binomial INGARCH (NB-INGARCH). Zhu (2012a) presented the Generalized Poisson INGARCH model, while Zhu (2012b) introduced the COM-Poisson INGARCH model. This paper also aims to propose a new model and hopes that the distribution is suitable for highly overdispersed time series data and features simple parameters.

To continue thinking from the perspective of Poisson distribution, among the two types of models that are suitable for integer-valued time series data, this paper tends to use the INGARCH model. Therefore, we focus on introducing the PINGARCH model proposed by Ferland (2006).  Fokianos (2009) considered nonlinear structures for PINGARCH model's mean process and utilized perturbation to establish the model's geometric ergodicity and the existence of arbitrary moments. 
 
Through the study and research of a series of models dealing with excessively integer-valued time series, the INGARCH model based on bell distribution is considered in this paper. Because Bell distribution is over-dispersed, bell distribution can handle some over-discrete data, and Bell distribution has only one parameter. It belongs to the single-parameter exponential family distribution, the parameter structure is simple, and it is convenient to put into use.

This paper is organized as follows. Section 2 introduces the bell distribution.
Section 3 provides the definition of the INGARCH model based on the Bell distribution and describes some statistical properties of the model. In Section 4,  the conditional maximum likelihood estimation (CML) method is employed to estimate the model's parameters. Simulation results are performed in section 5, and section 6 applies the model to two real examples.

\section{Bell Distribution}
The bell distribution with a single-parameter, is over-dispersed discrete distribution that belongs to the exponential family distribution. Its probability mass function is simple, and it is infinitely divisible.

To begin with is the introduction of the bell number, which was first proposed by BELL(1934) as below:
\begin{equation*}
 exp(e^{x}-1)=\sum_{n=0}^{\infty} \frac{B_n}{n!} x^n,x\in\mathbb{R},   
\end{equation*}
where $B_n$ is the bell number, which is given by:
\begin{equation}{\label{e2.1}}
    B_n=\frac{1}{e}\sum_{k=0}^{\infty}\frac{k^n}{k!}.
\end{equation}
It is obvious that the bell number $B_n$ is the n$th$ moments of a Poisson distribution with the parameter of 1.

The bell distribution is given below:\\
\textbf{Definition 1}: For a discrete random variable $Z$ with values in the set $\mathbb{N}
={0,1,2...}$, if its probability density function is as below:
\begin{equation}{\label{e2.2}}
  P_r(Z=z)=\frac{\theta^z e^{-e^\theta+1} B_z}{z!},z\in\mathbb{N}_0,   
\end{equation}
where $B_z$ is the bell number defined in (\ref{e2.1}). Then we call Z is a BELL distribution with parameter $\theta>0$, short for $Z \sim Bell(\theta)$.

If $Z \sim Bell(\theta)$, its probability generating function can be:
\begin{center}
    $G_Z(s)=E(s^{Z})=exp(e^{s \theta}-e^{\theta}), |s| \leq 1.$
\end{center}
The mean and variance of $Z$ are:
\begin{equation}{\label{e2.3}}
    E({Z})=\theta e^{\theta},Var({Z})=\theta (1+\theta) e^{\theta}.
\end{equation}
We can find that $Var(Z)/E(Z)=1+ \theta >1$, meaning the bell distribution is over-dispered. Thus it can account for some highly discrete data.

The properties of beLL distribution can be concluded as follows:
\begin{itemize}
\item Poisson distribution is not nested within the bell distribution family, but for smaller parameter values, the bell distribution approximates the Poisson distribution.
\item Bell distribution possesses identifiability, strong monotonicity, and infinite divisibility.
\item If the random variable $Z \sim Bell(\theta)$, ${Y_i, i=1,2...N}$ is the zero-truncated Poisson distribution with parameter $\theta$ and $Y_i \sim Poisson(e^{\theta}-1)$, then $Y_1+Y_2+...+Y_N$ has the same distribution with $Z$. 
\end{itemize}

\section{The proposed model}

 Suppose $\left\{X_{t}, {t} \geq 1\right\}$ denotes the  count time series, $\mathcal{F}_t$ is generated by $\sigma$-field $\left\{{X}_{{s}}, {s} \leq\right.$ $\left.t, \lambda_0\right\}, t \in \mathbb{Z}$, where the $\lambda_0$ denotes some initial values. Consider the below model:
\begin{equation}{\label{e3.1}}
 X_t \mid \mathcal{F}_{t-1} \sim \operatorname{Bell}\left(\lambda_t\right), \quad 
 \lambda_t=f\left(X_{t-1}, \ldots, X_{t-p}, \lambda_{t-1}, \ldots, \lambda_{t-q}\right),
\end{equation}
where $\alpha_0>0, \alpha_i \geq 0, \beta_j \geq 0, i=1, \ldots, p, j=1, \ldots, q, p \geq 1, q \geq 1$. Model (\ref{e3.1}) is denoted by BELL-INGARCH $(p, q)$. When $q=0$, the model is denoted by BELL-INGARCH$(p)$, and when $p=0$, it is denoted by BELL-INGARCH$(q)$. According to the properties of the bell distribution (\ref{e2.3}), we can obtain the approximate conditional expectation and variance can be obtained:
\begin{equation*}
{E}\left({X}_{{t}} \mid \mathcal{F}_{{t}-1}\right)= \lambda_{{t}} {e}^{\lambda_{{t}}}, \quad \operatorname{Var}\left({X}_{{t}} \mid \mathcal{F}_{{t}-1}\right) = \lambda_{{t}}\left(1+\lambda_{{t}}\right) {e}^{\lambda_{{t}}}
\end{equation*}
%Similar to Ferland et al.(2006), the BELL-INGARCH model was approximately stationary and possess finite first and second moments when $\sum_{{i}=1}^{{p}} \alpha_{{i}}+\sum_{{i}=1}^{{q}} \beta_{{j}}<1$.

When ${f}(\cdot)$ is a linear function,
\begin{equation}{\label{e3.2}}
\lambda_t=\alpha_0+\alpha_1 X_{t-1}+\beta_1 \lambda_{t-1}. 
\end{equation}

When ${f}(\cdot)$ is a non-linear function, 
\begin{equation}{\label{e3.3}}
\lambda_t=\frac{\alpha_0}{\left(1+\lambda_{t-1}\right)^\gamma}+\alpha_1 X_{t-1}+\beta_1 \lambda_{t-1}.  
\end{equation}
Suppose all the parameters $\alpha_0, \alpha_1, \beta_1, \gamma>0$, The introduction of $\gamma$ introduces non-linear perturbations. In this sense, the smaller the values of the $\gamma$, the closer equation (\ref{e3.3}) approximates equation (\ref{e3.2}), while larger parameter $\gamma$ values introduce stronger perturbations.

Assuming that model (\ref{e3.2}) is stationary, we can derive some important statistical properties of this model. For convenience, we define $Z_t$ is a conditional mean process with  $Z_t={E}\left({X}_{t}\mid \mathcal{F}_{{t}-1}\right)=\lambda_t e^{\lambda_t}$.
\begin{theorem}
    Suppose that $\left\{X_t, t \geq 1\right\}$ from an second-order stationary BELL-INGARCH $(p, q)$ model, $\sum_{i=1}^p \alpha_i+\sum_{i=1}^q \beta_j<1$, We can obtain:\\
(1) Unconditional expectation of $X_t$:
$$
\begin{aligned}
\mu & ={E}\left({X}_{{t}}\right)={E}\left[{E}\left({X}_{{t}} \mid \mathcal{F}_{{t}-1}\right)\right]={E}\left(\lambda_{{t}} {e}^{\lambda_{{t}}}\right)\\
& =\frac{\left(1-\sum_{j=1}^q \beta_j\right) \ln \alpha_0+\sum_{i=1}^p \alpha_i\left(\alpha_0-\ln \alpha_0\right)}{\sum_{i=1}^p \alpha_i} \cdot \exp \left(\frac{\left(1-\sum_{j=1}^q \beta_j\right) \ln \alpha_0+\sum_{i=1}^p \alpha_i\left(\alpha_0-\ln \alpha_0\right)}{\sum_{i=1}^p \alpha_i}\right);
\end{aligned}
$$
(2) Unconditional variance of $X_t$:
$$
\begin{gathered}
\sigma^2=\operatorname{Var}\left(X_t\right)=E\left(\operatorname{Var}\left(X_t \mid \mathcal{F}_{t-1}\right)\right)+\operatorname{Var}\left(E\left(X_t \mid \mathcal{F}_{t-1}\right)\right) \\
= E\left(\lambda_t\right)\left(1+E\left(\lambda_t\right)\right) \cdot e^{E\left(\lambda_t\right)}+\operatorname{Var}\left(\lambda_t\right) \cdot e^{\operatorname{Var}\left(\lambda_t\right)};
\end{gathered}
$$
(3) Covariance of $X_t$, and $Z_{t-k}$, $\operatorname{cov}\left(X_t, Z_{t-k}\right)$:
$$
\operatorname{cov}\left({X}_t, Z_{t-k}\right) =\left\{\begin{array}{ll}
\operatorname{cov}\left(\lambda_{{t}} {e}^{\lambda_{{t}}}, \lambda_{{t}-{k}} {e}^{\lambda_{t-{k}}}\right), & {k} \geq 0 \\
\operatorname{cov}\left({X}_{{t}}, {X}_{{t}-{k}}\right), & {k}<0
\end{array} ;\right.
$$
(4) Note that $\gamma_X(k)=\operatorname{cov}\left(X_t, X_{t-k}\right)$, then we can obtain $\operatorname{cov}\left(X_t, X_{t-k}\right)=\operatorname{cov}\left(\lambda_t {e}^{\lambda_t}, X_{t-k}\right)$. So, for $k \geq 1$,
$$
\begin{aligned}
\operatorname{cov}\left(Z_{{t}}, X_{t-k}\right) & =\operatorname{cov}\left(\alpha_0+\sum_{i=1}^p \alpha_i X_{t-i}+\sum_{j=1}^q \beta_j Z_{t-j}, X_{t-k}\right) \\
& = \sum_{i=1}^p \alpha_i \gamma_X(|k-i|)+\sum_{j=1}^{\min (k-1, q)} \beta_j \gamma_X(k-j)+\sum_{j=1}^q \beta_j \gamma_Z(j-k) .
\end{aligned}
$$
Then $\gamma_X(k)=\operatorname{cov}\left(X_t, X_{t-k}\right)=\operatorname{cov}\left(\lambda_{{t}} {e}^{\lambda_{{t}}}, X_{t-k}\right)$ can be obtained from the above.
\end{theorem}
\noindent{\textbf{Proof}}: The proofs for (1) and (2) are similar to Theorem 1 in Wei{\ss} (2009), and we will not reiterate them here. We will only provide the proof for (3), and (4) can be derived from (3).\\
Suppose $\mathcal{J}_t$ comes from the $\sigma-field$ generated from $\left\{\lambda_t, \lambda_{t-1}, \ldots\right\}$, Then we can get:
$$
{E}\left({X}_{{t}} \mid \mathcal{F}_{{t}-1}, \mathcal{J}_t\right)={E}\left({X}_{{t}} \mid \mathcal{F}_{{t}-1}\right)=\lambda_{{t}} {e}^{\lambda_{{t}}}.
$$
When ${k} \geq 0$, it is given as below:
$$
\begin{aligned}
\operatorname{cov}\left({X}_{{t}}-\lambda_{{t}} {e}^{\lambda_{{t}}}, \lambda_{{t}-{k}} {e}^{\lambda_{{t}-{k}}}\right) & ={E}\left[\left({X}_{{t}}-\lambda_{{t}} {e}^{\lambda_{{t}}}\right)\left(\lambda_{{t}-{k}} {e}^{\lambda_{t-{k}}}-\mu\right)\right] \\
& ={E}\left[\left(\lambda_{{t}-{k}} {e}^{\lambda_{{t}-{k}}}-\mu\right) {E}\left({X}_{{t}}-\lambda_{{t}} {e}^{\lambda_{{t}}} \mid \mathcal{J}_t\right)\right] \\
& ={E}\left[\left(\lambda_{{t}-{k}} {e}^{\lambda_{t-{k}}}-\mu\right)\left[{E}\left({E}\left({X}_{{t}} \mid \mathcal{F}_{{t}-1}, \mathcal{J}_t\right) \mid \mathcal{J}_t\right)-\lambda_{{t}} {e}^{\lambda_{{t}}}\right]\right] \\
& = {E}\left[\left(\lambda_{{t}-{k}} {e}^{\lambda_{t-{k}}}-\mu\right)\left[{E}\left(\lambda_{{t}} {e}^{\lambda_{{t}} \mid \mathcal{J}_t}\right)-\lambda_{{t}} {e}^{\lambda_{{t}}}\right]\right]=0,
\end{aligned}
$$
When ${k}<0$, it is given as below:
$$
\begin{aligned}
\operatorname{cov}\left({X}_{{t}}, {X}_{{t}-{k}}-\lambda_{{t}-{k}} {e}^{\lambda_{{t}-{k}}}\right) & ={E}\left[\left({X}_{{t}}-\mu\right)\left({X}_{{t}-{k}}-\lambda_{{t}-{k}} {e}^{\lambda_{{t}-{k}}}\right)\right] \\
& ={E}\left[\left({X}_{{t}}-\mu\right) {E}\left[\left({X}_{{t}-{k}}-\lambda_{{t}-{k}} {e}^{\lambda_{t-{k}}}\right) \mid \mathcal{F}_{{t}-{k}-1}\right]\right] \\
& = {E}\left[\left({X}_{{t}}-\mu\right)\left[\lambda_{{t}-{k}} {e}^{\lambda_{{t}-{k}}}-{E}\left(\lambda_{{t}-{k}} {e}^{\lambda_{{t}-{k}}} \mid \mathcal{F}_{{t}-{k}-1}\right)\right]\right]=0.
\end{aligned}
$$

\section{Estimation}

Suppose $\left\{{X}_{{t}}, {t} \geq 1\right\}$ denote the observed count time series satisfying the BELL-INGARCH $(p,q)$ model (\ref{e3.1}). $\theta=(\alpha_{0}, \alpha_{1}, \ldots, \alpha _{p}, \beta_{1}, \ldots, \beta_{q})^\intercal$ is the smallest value of average KL divergence $KL(\theta)$, $p^\ast=max(p,q)$, by $(\ref{e2.1})$ and $(\ref{e2.2})$ we can derive the probability density function of $X_t$:
\begin{center}
      $P(X_t=x_t\mid \lambda_t)=\frac{\lambda^{x_t}_{t} e^{(-e^{\lambda_{t}}+1)} B_{x_t}}{x_{t}!},$
\end{center}
where $B_{x_t}=\frac{1}{e} \sum_{k=0}^{\infty} \frac{k^{x_t}}{k !}$.

For the logarithmic likelihood function, where $n \in \mathbb{N}$ is the sample size:
$$
{l}(\theta)=\sum_{t=p^*}^{n}{l}_{t}(\theta)=\sum_{t=p^*}^{n}\left\{x_{t} \log \left(\lambda_{t}\right)+\left(-e^{\lambda_{t}}+1\right)+\log \left(B_{x_t}\right)-\log \left(x_{t} !\right)\right\} .
$$

The first derivative of ${l}_{t}(\theta)$ with respect to $\theta$ is given by:
$$
\frac{\partial {l}_{t}(\theta)}{\partial(\theta)}=\left[\frac{x_t}{\lambda_t}-\lambda_t e^{\lambda_t}\right] \frac{\partial \lambda_t}{\partial \theta}.
$$

The second derivative of ${l}_{{t}}(\theta)$ can be obtained as follows:
$$
\frac{\partial^2 l_t(\theta)}{\partial \theta \partial \theta^T} = \left[-\frac{x_t}{\lambda_t^2}-e^{\lambda_t}-\lambda_t^2 e^{\lambda_t}\right]\left(\frac{\partial \lambda_t}{\partial \theta}\right)\left(\frac{\partial \lambda_t}{\partial \theta}\right)^T + \left[\frac{x_t}{\lambda_t}-\lambda_t e^{\lambda_t}\right] \frac{\partial^2 \lambda_t}{\partial \theta \partial \theta^T}.
$$
Where $\partial\lambda_t / \partial\theta$ and $\partial^2 \lambda_t / \partial \theta \partial \theta^T$ are provided expressions in Fokianos et al. (2009). The CML estimate $\hat{\theta}$ is a solution to the score equation $\partial {l}(\theta) / \partial \theta =0$. Besides, the estimates of the model, denoted as $\hat{\theta}=(\widehat{\alpha}_0, \widehat{\alpha}_1, \ldots, \widehat{\alpha}_p, \widehat{\beta}_1, \ldots, \widehat{\beta}_q)^T$, can be obtained not only by solving the score equation but also by numerically maximizing the log-likelihood function. In this section, we primarily use the constrOptim function in the R language sofeware to directly numerically maximize the log-likelihood function, subject to the constraints $\alpha_0>0, \alpha_1 \geq 0, \beta_1 \geq 0, \alpha_1+\beta_1<1$. If we further differentiate the score equation, we can obtain the Hessian matrix of model (\ref{e3.1}):
$$
H_n(\theta)=-\sum_{t=p^*}^{n} \frac{\partial^2 l_{t}(\theta)}{\partial \theta \partial \theta^T}.
$$

\section{Simulations}
\subsection{BELL-INGARCH(1,1)}
\begin{table}[htbp]
	\label{T1}\caption{The estimated results of the BELL-INGARCH(1,1) model} \vspace{0.1in}
	\tabcolsep0.2in
	\begin{center}
		\small
		\begin{tabular}{cccccc}
			\hline  & ${n}$ & & $\widehat{a}_0$ & $\widehat{\alpha}_1$ & $\hat{\beta}_1$ \\
			\hline \multirow[t]{9}{*}{${Al}$} & 200 & Mean & 0.5806 & 0.0515 & 0.1433 \\
 & & MADE & 0.1173 & 0.0202 & 0.1506 \\
 & & MSE & 0.0242 & 0.0007 & 0.0450 \\
 & 500 & Mean & 0.5873 & 0.0558 & 0.1257 \\
 & & MADE & 0.1010 & 0.0134 & 0.1313 \\
 & & MSE & 0.0164 & 0.0003 & 0.0291 \\
 & 1000 & Mean & 0.6017 & 0.0586 & 0.0998 \\
 & & MADE & 0.0765 & 0.0083 & 0.0972\\
 & & MSE & 0.0086 & 0.0001 & 0.0142\\
\multirow[t]{9}{*}{${A} 2$} & 200 & Mean & 0.6244 & 0.0230 & 0.1769 \\
 & & MADE & 0.1437 & 0.0132 & 0.1793 \\
 & & MSE & 0.0412 & 0.0003 & 0.0667 \\
 & 500 & Mean & 0.6374 & 0.0246 & 0.1573 \\
 & & MADE & 0.1250 & 0.0087 & 0.1579 \\
 & & MSE & 0.0307 & 0.0001 & 0.0496 \\
 & 1000 & Mean & 0.6528 & 0.0241 & 0.1397 \\
 & & MADE & 0.1086 & 0.0059 & 0.1356 \\
 & & MSE & 0.0229 & 0.0001 & 0.0361 \\
\multirow[t]{9}{*}{${A} 3$} & 200 & Mean & 0.7703 & 0.0382 & 0.1233 \\
 & & MADE & 0.1192 & 0.0092 & 0.1239 \\
 & & MSE & 0.0271 & 0.0001 & 0.0308 \\
 & 500 & Mean & 0.7886 & 0.0390 & 0.1026 \\
 & & MADE & 0.0943 & 0.0066 & 0.0956 \\
 & & MSE & 0.0136 & 0.0001 & 0.0143 \\
 & 1000 & Mean & 0.8037 & 0.0395 & 0.0868 \\
 & & MADE & 0.0772 & 0.0044 & 0.0797 \\
 & & MSE & 0.0086 & 0.0000 & 0.0090 \\
\multirow[t]{9}{*}{${A} 4$} & 200 & Mean & 0.8835 & 0.0278 & 0.1408 \\
 & & MADE & 0.1581 & 0.0083 & 0.1432 \\
 & & MSE & 0.0384 & 0.0001 & 0.0326 \\
 & 500 & Mean & 0.9077 & 0.0293 & 0.1158 \\
 & & MADE & 0.1200 & 0.0051 & 0.1064 \\
 & & MSE & 0.0220 & 0.0000 & 0.0176 \\
 & 1000 & Mean & 0.9029 & 0.0298 & 0.1180 \\
 & & MADE & 0.1015 & 0.0036 & 0.0909 \\
 & & MSE & 0.0155 & 0.0000 & 0.0125 \\
\hline
		\end{tabular}
	\end{center}
\end{table}

In this section, numerical simulations based on the BELL-INGARCH(1,1) model are conducted to investigate the finite-sample properties of CML estimation. In terms of the selection of initial values, a random sampling method is employed from a uniform distribution centered around the true values with small errors. We select sample sizes of n = 200, 500, and 1000, and the simulation results are replicated by 500 times. We observe the mean, mean absolute deviation of estimates (MADE), and mean squared error (MSE) of the estimators for unknown parameters to assess the reasonableness and effectiveness of the estimates. Parameters under the following circumstances are considered:

The BELL - INGARCH(1,1) model when $\alpha_0>0, \alpha_1 \geq 0, \beta_1 \geq 0, \alpha_1+\beta_1<1$:

(A1) $\theta=(0.6,0.06,0.10)^{\top}$;

(A2) $\theta=(0.7,0.025,0.08)^{\top}$;

(A3) $\theta=(0.8,0.04,0.09)^\top$;

(A4) $\theta=(0.9,0.03,0.12)^{\top}$.

The estimated mean, mean absolute deviation of estimates, and mean squared error are given in Table 1. The Mean, represents the sum of all values in a dataset divided by the number of data points and reflects the central tendency of the data. MADE measures the average distance between the true values $\gamma_{{i}}$ and their estimated values $\hat{\gamma}_i$, providing a good reflection of the actual prediction errors. On the other hand, MSE calculates the average of the squared differences between the true values $\gamma_{{i}}$ and their estimated values $\hat{y}_i$, indicating the degree of data variation. A smaller MSE value suggests a higher level of precision in describing experimental data. The formulas MADE and MSE are provided below:

$$
\text { MADE }=\frac{1}{{n}} \sum_{{i}=1}^{{n}}\left|\hat{\gamma}_i-\gamma_{{i}}\right|, \quad \text { MSE }=\frac{1}{{n}} \sum_{{i}=1}^{{n}}\left(\hat{\gamma}_i-\gamma_{{i}}\right)^2.
$$

Table 1 reveals that the results of the CML estimation are reasonable. With an increasing sample size, the estimated mean from the CML method approaches the true parameter value, and the corresponding MADE and MSE become progressively smaller. This indicates an improvement in estimation accuracy. Furthermore, the estimated quantities exhibit small MADE and MSE values, highlighting the effectiveness of CML estimation.
\subsection{NBELL-GARCH(1,1)}
Further exploring the finite sample properties of CML (Conditional Maximum Likelihood) estimation through numerical simulations based on the NBELL-INGARCH(1,1) model. The initial values are still randomly selected from a uniform distribution centered around the true values with small errors. The simulation is repeated 500 times, considering sample sizes of $\mathrm{n}=200$, 500, and 1000. We observe the reasonability and effectiveness of the estimates through the mean, MAE (Mean Absolute Deviation), and MSE (Mean Squared Error) of the unknown parameter estimates. The parameters under consideration are $\alpha_0, \gamma > 0, \alpha_1, \beta_1 \geq 0, \alpha_1 + \beta_1 < 1$, and $\theta=(\alpha_0, \alpha_1, \beta_1, \gamma)^T$:

(B1) $\theta=(0.4,0.04,0.1,0.5)^{\top}$;

(B2) $\theta=(0.5,0.05,0.2,0.8)^{\top}$;

(B3) $\theta=(0.7,0.06,0.3,1.0)^{\top}$;

(B4) $\theta=(0.4,0.07,0.4,2.0)^{\top}$.

From Table 2, it is evident that under the numerical simulation of the NBELL-INGARCH $(1,1)$ model, CML estimates remain reasonable. With an increase in sample size, the average values estimated by the CML method converge to the true parameter values, and the corresponding MADE and MSE decrease, indicating improved estimation accuracy. The estimates also exhibit small MADE and MSE, affirming the effectiveness of CML estimation.
\begin{table}[htbp]
	\label{T1}\caption{The estimated results of the NBELL-INGARCH(1,1) model} \vspace{0.1in}
	\tabcolsep0.2in
	\begin{center}
		\small
		\begin{tabular}{ccccccc}
			\hline  & ${n}$ & & $\widehat{a}_0$ & $\widehat{\alpha}_1$ & $\hat{\beta}_1$ & $\hat{\gamma}$\\
			\hline \multirow[t]{9}{*}{${Bl}$} & 200 & Mean & 0.4084 & 0.0490 & 0.1049 & 0.5089 \\
			& & MADE & 0.0142 & 0.0116 & 0.0120 & 0.0150 \\
			& & MSE & 0.0005 & 0.0003 & 0.0005 & 0.0005 \\
			& 500 & Mean & 0.4063 & 0.0469 & 0.0993 & 0.5061 \\
			& & MADE & 0.0092 & 0.0094 & 0.0085 & 0.0094 \\
			& & MSE & 0.0002 & 0.0002 & 0.0002 & 0.0002 \\
			& 1000 & Mean & 0.4058 & 0.0448 & 0.1005 & 0.5055 \\
			& & MADE & 0.0083 & 0.0073 & 0.0066 & 0.0089 \\
			& & MSE & 0.0001 & 0.0001 & 0.0001 & 0.0002 \\
			\multirow[t]{9}{*}{${B} 2$} & 200 & Mean & 0.5116 & 0.0589 & 0.1954 & 0.8119 \\
			& & MADE & 0.0145 & 0.0116 & 0.0144 & 0.0151\\
			& & MSE & 0.0006 & 0.0002 & 0.0009 & 0.0006 \\
			& 500 & Mean & 0.5068 & 0.0569 & 0.1943 & 0.8070 \\
			& & MADE & 0.0096 & 0.0090 & 0.0106 & 0.0099 \\
			& & MSE & 0.0002 & 0.0001 & 0.0004 & 0.0003 \\
			& 1000 & Mean & 0.5060 & 0.0556 & 0.1959 & 0.8064 \\
			& & MADE & 0.0089 & 0.0075 & 0.0094 & 0.0093 \\
			& & MSE & 0.0002 & 0.0001 & 0.0003 & 0.0002  \\
			\multirow[t]{9}{*}{${B} 3$} & 200 & Mean & 0.7062 & 0.0655 & 0.2982 & 1.0066 \\
			& & MADE & 0.0096 & 0.0074 & 0.0092 & 0.0097 \\
			& & MSE & 0.0006 & 0.0001 & 0.0005 & 0.0006 \\
			& 500 & Mean & 0.7042 & 0.0639 & 0.2992 & 1.0035 \\
			& & MADE & 0.0071 & 0.0062 & 0.0068 & 0.0074\\
			& & MSE & 0.0002 & 0.0001 & 0.0002 & 0.0002  \\
			& 1000 & Mean & 0.7029 & 0.0643 & 0.2989 & 1.0027 \\
			& & MADE & 0.0058 & 0.0054 & 0.0062 & 0.0062 \\
			& & MSE & 0.0001 & 0.0000 & 0.0001 & 0.0001 \\
			\multirow[t]{9}{*}{${B} 4$} & 200 & Mean & 0.4055 & 0.0729 & 0.3523 & 2.0052 \\
			& & MADE & 0.0223 & 0.0138 & 0.0491 & 0.0222 \\
			& & MSE & 0.0018 & 0.0004 & 0.0055 & 0.0018 \\
			& 500 & Mean & 0.4030 & 0.0729 & 0.3694 & 2.0034 \\
			& & MADE & 0.0112 & 0.0085 & 0.0328 & 0.0121  \\
			& & MSE & 0.0007 & 0.0001 & 0.0027 & 0.0007 \\
			& 1000 & Mean & 0.4031 & 0.0721 & 0.3746 & 2.0027 \\
			& & MADE & 0.0080 & 0.0068 & 0.0271 & 0.0080 \\
			& & MSE & 0.0001 & 0.0001 & 0.0017 & 0.0001 \\
			\hline
		\end{tabular}
	\end{center}
\end{table}

\section{Empirical Applications}

In this section, the BELL-INGARCH $(1,1)$ model is employed to fit two distinct types of empirical datasets. A comparative analysis of the goodness of fit is conducted by contrasting the performance of this model against the PINGARCH $(1,1)$ model proposed by Ferland et al. (2006) and the NB-INGARCH $(1,1)$ model proposed by Christou and Fokianos (2014) for these two classes of datasets. For the INGARCH model, the constrOptim function in the R programming language is utilized to maximize the log-likelihood function. The imposed constraints for this optimization are defined as $\alpha_0>0, \alpha_1 \geq 0, \beta_1 \geq 0, \alpha_1+\beta_1<1$.

Next, we introduce the NB-INGARCH model. The NB-INGARCH model in time series is a non-linear equation-based time-varying volatility model used to characterize the volatility of time series data, particularly employed in the financial domain for modeling non-negative data such as stock returns or trading volumes. This model takes into account the discreteness and volatility of the data, allowing for the capture of heteroscedasticity within the dataset.

\subsection{Daily Download Count of TeX Editor}

\begin{figure}[htbp]
	\centering
	\includegraphics[width=0.7\textwidth]{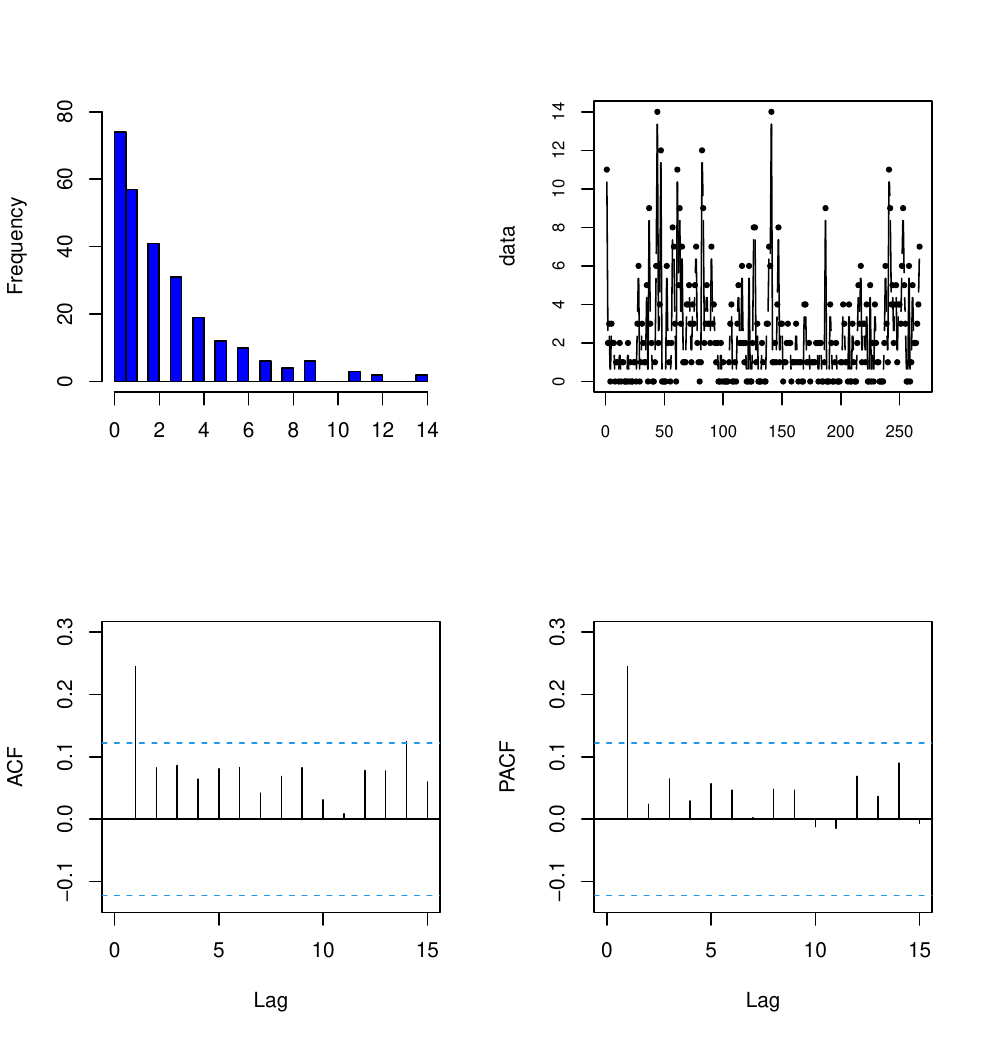}
	\caption{Histogram, Path Diagram, ACF, and PACF of Daily Downloads for TEX Editor} 
	\label{fig:figure1} 
\end{figure}

\begin{table}[htbp]
	\centering
	\caption{Parameter Estimation, Log-Likelihood Function, AIC, and BIC for Daily Download Counts of TEX Editor}
	\begin{tabular}{ccccc}
		&&&&\\\hline 
		Model & Parameters & Log-Likelihood & AIC & BIC \\
		\hline & $\widehat{\alpha}_0=1.3942$ & & & \\
		PINGARCH $(1,1)$ & $\widehat{\alpha}_1=0.1369$ & -623.0368 & 1252.0737 & 1262.8354 \\
		& $\widehat{\beta}_1=0.2730$ & & & \\
		\hline & $\widehat{\alpha}_0=1.3942$ & & & \\
		& $\widehat{\alpha}_1=0.1369$ & & & \\
		NB-INGARCH $(1,1)$ & $\widehat{\beta}_1=0.2730$ & -623.0368 & 1256.0737 & 1274.0099 \\
		& $\widehat{v}_1=1.2340$ & & & \\
		& $\hat{v}_2=1.2079$ & & & \\
		\hline & $\widehat{\alpha_0}=0.4480$ & & & \\
		BELL-INGARCH $(1,1)$ & $\widehat{\alpha_1}=0.0406$ & -544.4992 & 1094.9984 & 1105.7602 \\
		& $\widehat{\beta_1}=0.4079$ & & & \\
		\hline
	\end{tabular}
\end{table}

\begin{figure}[htbp]
	\centering
	\includegraphics[width=0.7\textwidth]{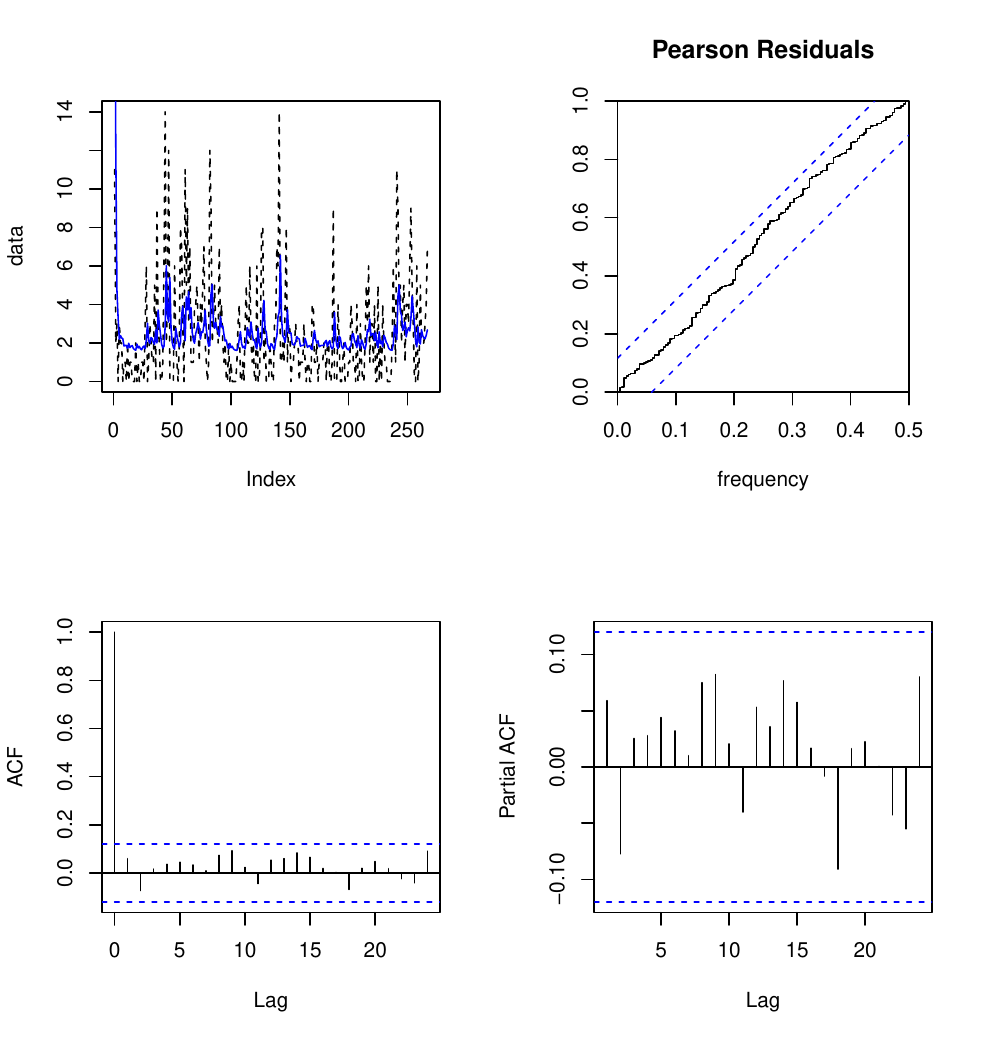}
	\caption{The observed and predicted values of daily downloads for the TEX editor, the residuals fitted with the BELLINGARCH model, and the ACF and PACF of residuals.} 
	\label{fig:figure2} 
\end{figure} 

In the first example, the dataset considered is proposed by Weiss (2008). This dataset comprises integer-valued time series data, representing the daily download counts of the TEX editor from June 2006 to February 2007 (\(T=267\)). The length of the analyzed data is 267, with mean and variance values of 2.4007 and 7.5343, respectively. Given that the variance exceeds the mean, indicating overdispersion, it is reasonable to hypothesize that the BELL-INGARCH \((1,1)\) model is suitable for fitting this dataset.

As illustrated in Figure 1, the figure presents the histogram, path diagram, ACF plot, and PACF plot for the daily download counts dataset of the TEX editor. The dataset is fitted using the PINGARCH $(1,1)$ model, NB2-INGARCH $(1,1)$ model, and BELL-INGARCH $(1,1)$ model. Table 2 summarizes the parameter estimates, log-likelihood function, AIC, and BIC for the daily download counts of the TEX editor resulting from the fitting process.

The Akaike Information Criterion (AIC) is a statistical model evaluation method. It can be expressed as \( {AIC} = 2k - 2 \ln(L) \), where \( k \) is the number of parameters, and \( L \) is the likelihood function. Built upon the concept of entropy, AIC balances the complexity of the estimated model against its goodness-of-fit to the data. The Bayesian Information Criterion (BIC) is also a statistical model evaluation method with a larger penalty term than AIC. Specifically, considering the sample size \( n \), BIC can be expressed as \( {BIC} = k \ln(n) - 2 \ln(L) \). When \( n \) is large, \( k \ln(n) \geq 2k \), so BIC effectively prevents excessive model complexity caused by excessively high model accuracy when there is a large sample size. Based on the AIC and BIC criteria, smaller AIC and BIC values indicate a better-fitting model. Thus, from Table 6-1, it can be concluded that, compared to other models, the BELL-INGARCH $(1,1)$ model demonstrates a better fit.

For the above model, Pearson residual analysis was also considered. The standard definition of Pearson residual is given by:
\[
e_t=\frac{X_t-{E}\left({Y}_{{t}} \mid \mathcal{F}_{{t}-1}\right)}{\sqrt{\operatorname{Var}\left({Y}_{{t}} \mid \mathcal{F}_{{t}-1}\right)}}.
\]

Upon calculation, the mean and standard error of Pearson residuals for the BELL-INGARCH $(1,1)$ model are 0.0170 and 1.2254, respectively. In Figure 6-2, the first plot depicts the observed values of daily downloads for the TEX editor as a black dashed line, while the blue solid line represents the predicted values defined as $\hat{Y}_{{t}}=\hat{\lambda}_{{t}} {e}^{\hat{\lambda}_{{t}}}$. From the figure, it is evident that the defined predictions effectively capture the upward and downward trends in the observed process. The second plot in Figure 6 displays the cumulative periodogram of Pearson residuals for the fitted model of weekly sales volume of a certain product. The dashed line represents the Kolmogorov-Smirnov boundary with a level of $\alpha=0.05$. It can be observed from the plot that the model fits the dataset well. The third and fourth plots in Figure 2 depict the autocorrelation function (ACF) and partial autocorrelation function (PACF) of Pearson residuals. The ACF and PACF of residuals suggest that the residuals are uncorrelated, further indicating the model's satisfactory fitting performance.

\subsection{Weekly Sales Volume of a Certain Soap Product}

In the second example, the BELL-INGARCH(1,1) model is employed to analyze the weekly sales volume of a particular fatty product in a supermarket from 1989 to 1994. The data is sourced from the \href{https://www.chicagobooth.edu/research/kilts/datasets/dominicks}{Kilts Marketing Center database} at the University of Chicago Booth School of Business. The specific product under consideration is "Zest White Water 15 oz.," with product code 3700031165, and store number 67. The sequence has a length of 242, with a mean and variance of 5.4421 and 15.4012, respectively. Given that the variance exceeds the mean, indicating pronounced volatility in the dataset, the BELL-INGARCH(1,1) model is deemed suitable for fitting this data. Figure 3 illustrates the histogram, path diagram, ACF, and PACF of the weekly sales volume for the soap product.

Using the PINGARCH $(1,1)$ model, NB2-INGARCH $(1,1)$ model, and BELL$\operatorname{INGARCH}(1,1)$ model to fit this dataset, we compute the parameter estimates, log-likelihood function, {AIC}, and {BIC}  for the weekly sales volume of the fertilized products, as shown in Table 3. Based on the {AIC} and {BIC} criteria, it can be inferred from Table 3 that the BELL-INGARCH $(1,1)$ model exhibits the best fitting performance.

This analysis also takes into account Pearson residual analysis, the mean and standard deviation of Pearson residuals for the BELL-INGARCH $(1,1)$ model are -0.0674 and 1.3198, respectively. From the first plot in Figure 4, it is evident that the model successfully captures the trend of fluctuations in the dataset. The second plot in Figure 4 illustrates the effective fitting of the model to the dataset. The third and fourth plots in Figure 4 reveal that the residuals are uncorrelated, providing further evidence of the model's strong fit to the dataset.

\begin{figure}[htbp]
	\centering
	\includegraphics[width=0.7\textwidth]{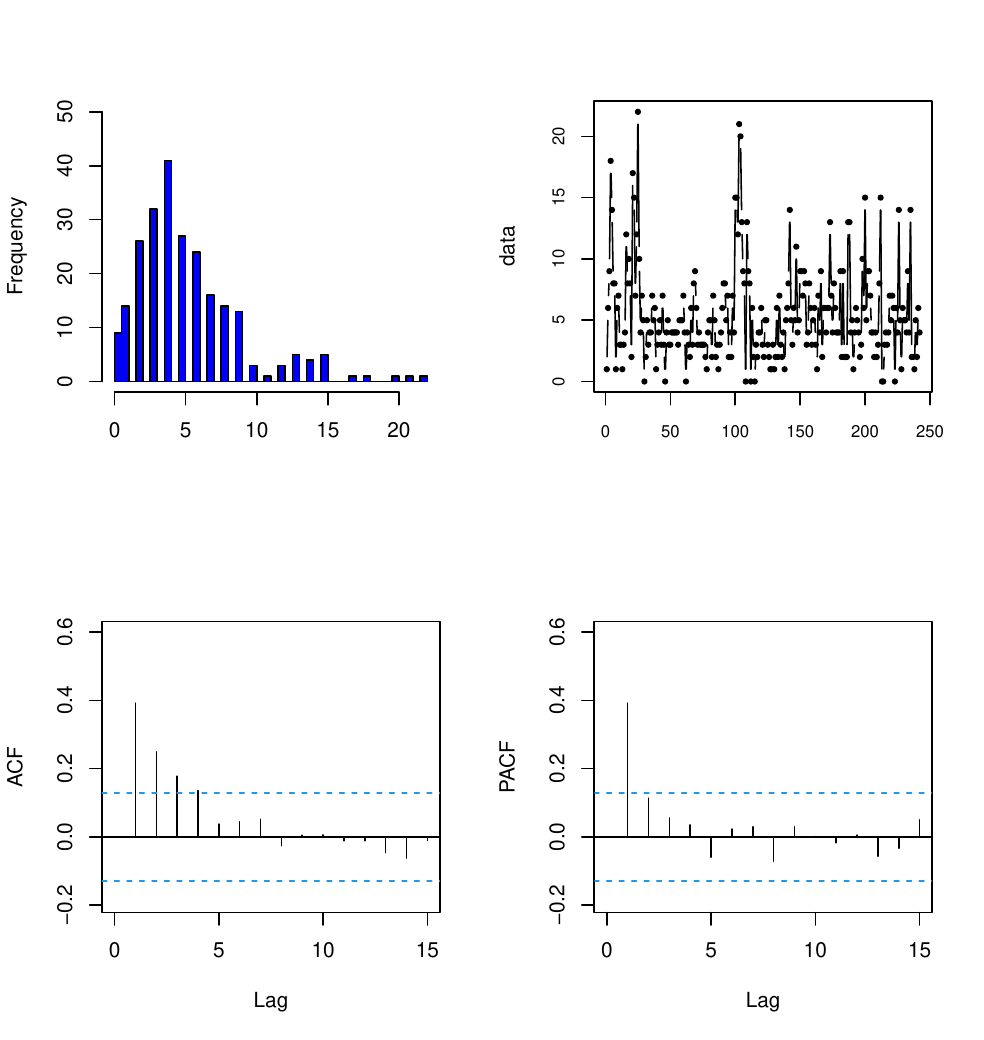}
	\caption{Histogram, Path Diagram, ACF, and PACF of Weekly Sales Volume of a Certain Soap Product} 
	\label{fig:figure3} 
\end{figure}

\begin{table}[htbp]
	\centering
	\caption{Parameter Estimation, Log-Likelihood Function, AIC, and BIC for Weekly Sales Volume of a Certain Soap Product}
	\begin{tabular}{ccccc}
		&&&&\\\hline 
		Model & Parameters & Log-Likelihood & AIC & BIC \\
		\hline & $\widehat{\alpha}_0=1.8055$ & & & \\
		PINGARCH $(1,1)$ & $\widehat{\alpha}_1=0.3637$ & -660.0511 & 1326.1023 & 1336.5691 \\
		& $\widehat{\beta}_1=0.3053$ & & & \\
		\hline & $\widehat{\alpha}_0=1.8055$ & & & \\
		& $\widehat{\alpha}_1=0.3637$ & & & \\
		NB-INGARCH $(1,1)$ & $\widehat{\beta}_1=0.3053$ & -660.0511 & 1330.1023 & 1347.5469 \\
		& $\widehat{v}_1=4.5011$ & & & \\
		& $\hat{v}_2=4.3488$ & & & \\
		\hline & $\widehat{\alpha_0}=0.9928$ & & & \\
		BELL-INGARCH $(1,1)$ & $\widehat{\alpha_1}=0.0335$ & -615.7368 & 1237.4735 & 1247.9403 \\
		& $\widehat{\beta_1}=0.1374$ & & & \\
		\hline
	\end{tabular}
\end{table}

\begin{figure}[htbp]
	\centering
	\includegraphics[width=0.7\textwidth]{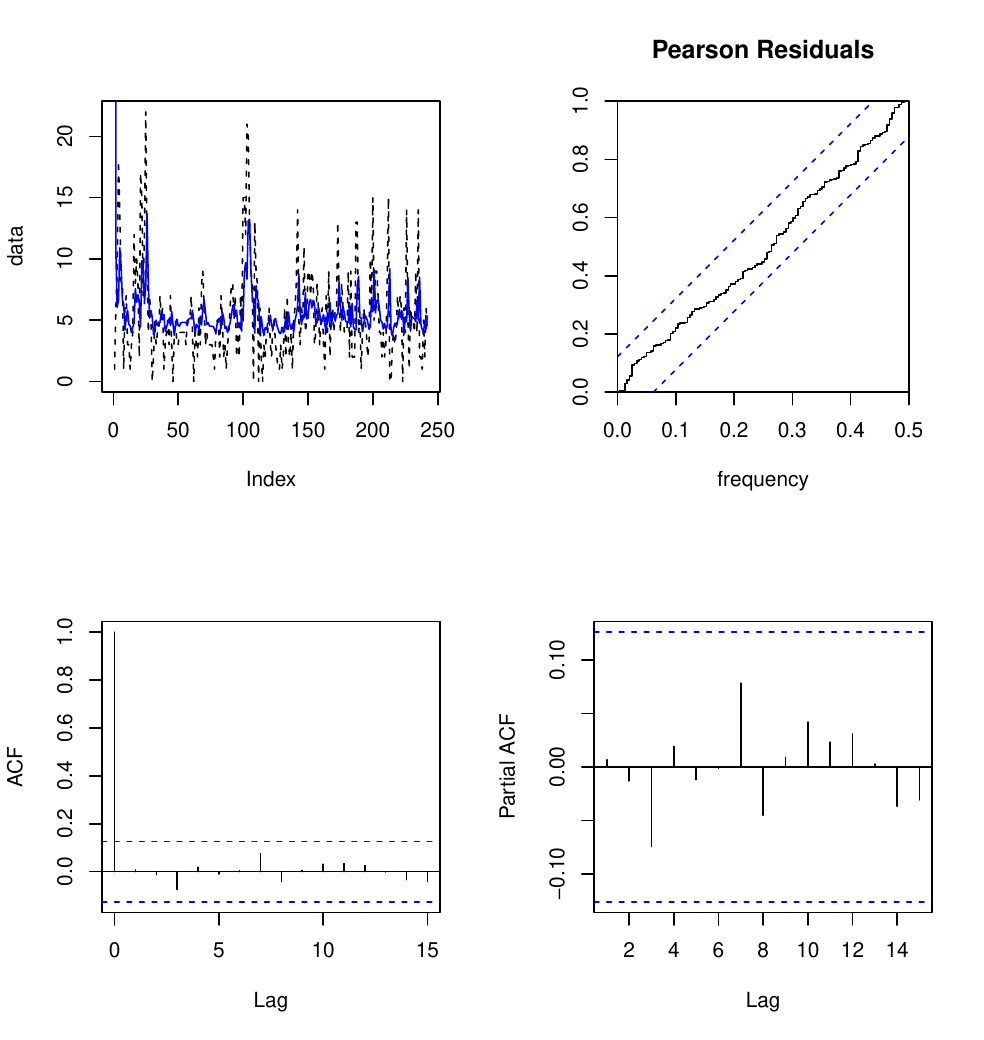}
	\caption{The observed and predicted values of Weekly Sales Volume of a Certain Soap Product, the cumulative cycle plot of Pearson residuals fitted with the BELLINGARCH $(1,1)$ model for Weekly Sales Volume of a Certain Soap Product, and ACF and PACF of Pearson residuals.} 
	\label{fig:figure4} 
\end{figure}

\section{Conclusion}
This paper proposed the new INGARCH model based on Bell distribution. The conditional maximum likelihood  estimation method is used to obtain the estimators of parameters. Numerical simulations confirm the finite sample properties of the estimation of unknown parameters. Finally, the  model is applied in the two real examples. Compared with the existing models, the proposed model is more simple and applicable.

\vspace{3mm}
{\Large\bf Acknowledgements}
\vspace{1mm}

This work is supported by the Project Funded by the Priority Academic Program Development of Jiangsu Higher Education Institutions and the Teacher's Research Support Project Foundation of Jiangsu Normal University (No. 21XFRS022).

\end{document}